\documentclass[11pt]{article}
\usepackage{amsmath, epsfig, cite, amssymb, amsfonts,latexsym}
\newcommand{\qed}{{\hfill\rule{4pt}{7pt}}}

\def\MX{\mathbb{X}}
\def\MY{\mathbb{Y}}
\def\MZ{\mathbb{Z}}

\newtheorem{thm}{Theorem}[section]

\newtheorem{lem}[thm]{Lemma}

\def\pf{\noindent {\it Proof.} }

\makeatletter \@addtoreset{equation}{section} \makeatother

\begin{document}
\begin{center}
{\Large \bf  $q$-Identities related to overpartitions and divisor
functions}
\end{center}

 \vspace{0.3cm}
 \begin{center}
   {\bf Amy M. Fu}\\
   Center for Combinatorics, LPMC\\
   Nankai University, Tianjin 300071, P.R. China\\
   Email: fmu@eyou.com
  \end{center}
 \begin{center}
  {\bf Alain Lascoux}\\
  Nankai University, Tianjin 300071, P.R. China\\
  Email: Alain.Lascoux@univ-mlv.fr\\
  CNRS, IGM Universit\'e  de
  Marne-la-Vall\'ee\\
  77454 Marne-la-Vall\'ee Cedex, France\\
  \end{center}

  \vspace{1cm}

\noindent{\bf Abstract.} We prove some $q$-Identities 
related to overpartitions and divisor functions

\section{Introduction}

In this note,
using the classical notations $(z;q)_i=(1-z)\cdots(1-zq^{i-1})$,
and $\displaystyle{{n \brack i}}=(q;q)_n/(q;q)_i(q;q)_{n-i}$,
 we prove the following two identities~:

\begin{thm}   For any pair positive integers $m, n$, one has 
\begin{multline}\label{EQ1}
\sum_{i=1}^n{n \brack
i}\frac{(-1)^{i-1}(x+1)\ldots(x+q^{i-1})}{(1-q^i)^m}q^{mi}
\\=\sum_{i=1}^{n}\frac{(-1)^{i-1}(x^i-(-1)^i)}{1-q^i}q^i
 \sum_{i\leq i_2\leq\ldots\leq i_m\leq
n}\frac{q^{\sum_{j=2}^mi_j}}{\prod_{j=2}^m(1-q^{i_j})},
\end{multline}

\begin{equation}\label{EQ2}
\frac{(z;q)_{n+1}}{(q;q)_n}\sum_{i=0}^n {n \brack
i}\frac{(-1)^{i-1}(x+1)\ldots(x+q^{i-1})}{1-zq^i}q^i
=\sum_{i=0}^n(-1)^{i-1}\frac{(z;q)_i}{(q;q)_i}x^iq^i \, .
\end{equation}
\end{thm}

In the next section, we shall show that (\ref{EQ1}) and
(\ref{EQ2}) can be obtained from the Newton interpolation in
points $\{-1,-q,-q^2,\ldots\}$, using complete functions  in the
variables $\{q/(1-q), q^2/(1-q^2), \cdots\}$.

Using the same methods, we have already given in \cite{FuLa}
an identity which generalizes the case $x=0$ of (\ref{EQ1}).  

Given $\MX=\{x_1, x_2, \ldots\}$, Newton gave the following
interpolation formula, for any function $f(x)$:
\begin{eqnarray*}
f(x)=f(x_1)+f\partial_1(x-x_1)+f\partial_1\partial_2(x-x_1)(x-x_2)+\ldots,
\end{eqnarray*}
where $\partial_i$, acting on its left, is defined by
$$
f(x_1,\ldots,x_i,x_{i+1},\ldots)\partial_i=\frac{f(\ldots,x_i,x_{i+1},\ldots)-f(\ldots,x_{i+1},x_i,\ldots)}{x_i-x_{i+1}}.
$$

Taking in particular  $f(x)=x^n$, we have
\begin{equation}\label{Com}
x_1^n \partial_1\ldots \partial_i=h_{n-i}(x_1,x_2,\ldots,x_{i+1}),
\end{equation}
where $h_k$ is the complete function of degree $k$, defined by
$$
h_{k}(x_1,x_2,\ldots, x_n)=\sum_{1 \leq i_1 \leq i_2 \leq \ldots
\leq i_k \leq n}x_{i_1}x_{i_2}\ldots x_{i_k}.
$$
Recall the following properties of $h_k$:
\begin{itemize}
\item[1.] The generating function of $h_k$
is:
\begin{equation}\label{Gene}
\sum_{k=0}^{\infty}h_{k}(x_1,\ldots,x_n)t^k=\frac{1}{(1-tx_1)\ldots(1-tx_n)}.
\end{equation}

 \item[2.] More generally, given two alphabets $\MX$ and $\MY$, 
then the generating functions of $h_k(\MX+\MY)$
and $h_k(\MX-\MY)$ are:
\begin{eqnarray}
\sum_{k=0}^{\infty}h_k(\MX+\MY)t^k &=& \frac{1}{\prod_{x \in
\MX}(1-xt)\prod_{y \in \MY}(1-yt)}, \nonumber\\
\sum_{k=0}^{\infty}h_k(\MX-\MY)t^k &=& \frac{\prod_{y \in
\MY}(1-yt)}{\prod_{x \in \MX}(1-xt)}\label{gened}.
\end{eqnarray}

As a consequence, one has:
\begin{equation}\label{expan}
h_n(\MX + \MY)=\sum_{k=0}^n h_k(\MX)h_{n-k}(\MY).
\end{equation}

 \item[3.] Given $\{x_1,x_2,\ldots,x_n\}$, and a positive integer $m$, we
 have:
\begin{equation}\label{Pas}
\sum_{i=1}^n x_i h_{m-1}(x_i, x_{i+1},\ldots,
x_{n})=h_m(x_1,x_2,\ldots,x_n).
\end{equation}
\end{itemize}

Taking  $\MX=\{-1,-q,-q^2, \ldots\}$, it is easy to check from
(\ref{Com}) and (\ref{Gene}):
$$
x_1^n\partial_1\ldots\partial_i=h_{n-i}(-1,-q,\ldots,-q^{i})=(-1)^{n-i}{n
\brack i}.
$$

\section{Proofs of (\ref{EQ1}) and (\ref{EQ2})}

 The Gauss polynomials $\displaystyle{{n \brack k}}$ satisfy the following recursion
 (cf.\cite{Andrews1}):
\begin{equation}\label{Pascal2}
\sum_{j=0}^n{m+j \brack m}q^j={n+m+1 \brack m+1}.
\end{equation}

In this paper, we need the following more general relations.

\begin{lem} \label{LEM}
Let $k$, $m$ and $n$ be nonnegative integers. Then we have the
following formulas:
\begin{equation}\label{Pascal1}
\sum_{i=k}^n{i \brack k}\frac{q^i}{1-q^i}\sum_{i \leq i_2 \ldots
\leq i_m \leq n}\frac{q^{\sum_{j=2}^m
i_j}}{\prod_{j=2}^m(1-q^{i_j})}={n \brack k}
\frac{q^{km}}{(1-q^k)^m},
\end{equation}
and
\begin{equation}\label{Pasz}
\sum_{i=0}^n\frac{(z;q)_i}{(q;q)_i}q^i=\frac{(zq;q)_n}{(q;q)_n}.
\end{equation}
\end{lem}

\pf\  Taking $\MX=\{1,q,\ldots, q^l\}$ and
$\MY=\{q^{l+1},q^{l+2},\ldots\}$, we obtain from (\ref{gened})
 and (\ref{expan}):
\begin{equation}\label{Aqq1}
\frac{1}{(q;q)_n}=h_n(\MX+\MY)=\sum_{i=0}^n
h_i(\MY)h_{n-i}(\MX)=\sum_{i=0}^n\frac{1}{(q;q)_i}{n-i+l \brack
l}q^{(l+1)i}.
\end{equation}

Letting $f(m)$ be the left side of (\ref{Pascal1}), then we have:
\begin{eqnarray*}
\sum_{m=1}^{\infty}f(m)z^m&=&\sum_{m=1}^{\infty}\sum_{i=k}^n{i
\brack k}\frac{q^i}{1-q^i}h_{m-1}\left(\frac{q^i}{1-q^i},
\frac{q^{i+1}}{1-q^{i+1}},\ldots, \frac{q^n}{1-q^n}\right)z^m\\
&\overset{(\ref{Gene})}{=}&z\sum_{i=k}^n{i \brack
k}\frac{q^i}{1-q^i}\frac{1}{(1-q^iz/(1-q^i))\ldots(1-q^nz/(1-q^n))}\\
&=&z\sum_{i=k}^n{i \brack
k}\frac{q^i}{1-q^i}\sum_{l=0}^{\infty}(q^i;q)_{n-i+1}{n-i+l \brack l}(q^i(1+z))^l\\
&=&z{n \brack
k}\sum_{i=k}^n\frac{(q;q)_{n-k}}{(q;q)_{i-k}}\sum_{l=0}^{\infty}{n-i+l
\brack l}q^{(l+1)i}\sum_{m=0}^l\binom{l}{m}z^m  \\   
&=&z{n \brack
k}\sum_{l=0}^{\infty}\sum_{m=0}^l\binom{l}{m}z^m\sum_{i=k}^n\frac{(q;q)_{n-k}}{(q;q)_{i-k}}{n-i+l
\brack l}q^{(l+1)i}\\
&\overset{(\ref{Aqq1})}{=}&{n \brack
k}\sum_{l=0}^{\infty}\sum_{m=0}^l\binom{l}{m}z^{m+1}q^{(l+1)k} 
  \end{eqnarray*}   
$$  = {n \brack k}\frac{q^{k}z}{1-q^k(1+z)}\\
=\sum_{m=1}^{\infty}{n \brack k} \frac{q^{km}}{(1-q^k)^m}z^m \, . $$

Taking $\MX=\{1\}$, $\MY=\{q,q^2,\ldots\}$ and $\MZ=\{zq, zq^2,
\ldots\}$, we get from (\ref{gened}) and (\ref{expan}), the proof
of (\ref{Pasz}):
\begin{multline*}
\frac{(zq;q)_n}{(q;q)_n}=h_{n}((\MX+\MY)-\MZ)\\=h_{n}(\MX+(\MY-\MZ))
=\sum_{i=0}^n
h_{i}(\MY-\MZ)h_{n-i}(\MX)=\sum_{i=0}^n\frac{(z;q)_i}{(q;q)_i}q^i.
\end{multline*}

\qed

Taking
\begin{eqnarray*}
f(x)&=&\sum_{i=1}^{n}\frac{(-1)^{i-1}(x^i-(-1)^i)}{1-q^i}q^i
\sum_{i\leq i_2\leq\ldots\leq i_m\leq
n}\frac{q^{\sum_{j=2}^mi_j}}{\prod_{j=2}^m(1-q^{i_j})},
\end{eqnarray*}
and $ \MX=\{-1,-q,-q^2,\ldots\}, $
then we have,
\begin{eqnarray*}
f(x)&=&f(x_1)+\sum_{k=1}^n f(x_1)\partial_1\ldots\partial_k
\, (x+1)\ldots(x+q^{k-1})\\
&=&\sum_{k=1}^n\sum_{i=k}^n(-1)^{k-1}{i \brack
k}\frac{q^i}{1-q^i}\sum_{i\leq i_2\leq\ldots\leq i_m\leq
n}\frac{q^{\sum_{j=2}^mi_j}}{\prod_{j=2}^m(1-q^{i_j})}(x+1)\ldots(x+q^{k-1})\\
&=&\sum_{k=1}^n(-1)^{k-1}(x+1)\ldots(x+q^{k-1})\sum_{i=k}^n{i
\brack k}\frac{q^i}{1-q^i}\sum_{i\leq i_2\leq\ldots\leq i_m\leq
n}\frac{q^{\sum_{j=2}^mi_j}}{\prod_{j=2}^m(1-q^{i_j})}\\
&=&\sum_{k=1}^n{n \brack
k}\frac{(-1)^{k-1}(x+1)\ldots(x+q^{k-1})}{(1-q^k)^m}q^{mk},
\end{eqnarray*}
as stated in (\ref{EQ1}). \qed

Taking
\begin{eqnarray*}
f(x)=\sum_{i=0}^{n}(-1)^{i-1}\frac{(z;q)_i}{(q;q)_i}x^{i}q^i,
\quad  \text{and} \quad  \MX=\{-1,-q,-q^2, \ldots\},
\end{eqnarray*}
we have,
\begin{eqnarray*}
f(x)&=&\sum_{k=0}^nf(x_1)\partial_1\ldots\partial_k \, (x+1)\ldots(x+q^{k-1})\\
&=&\sum_{k=0}^n (x+1)\ldots(x+q^{k-1}) \sum_{i=k}^n(-1)^{k-1}{i
\brack k}\frac{(z;q)_i}{(q;q)_i}q^i\\
&=&\sum_{k=0}^n(-1)^{k-1}(x+1)\ldots(x+q^{k-1})q^k\frac{(z;q)_k}{(q;q)_k}
\sum_{i=k}^n\frac{(zq^k;q)_{i-k}}{(q;q)_{i-k}}q^{i-k} 
\end{eqnarray*}
\begin{eqnarray*}
\phantom{f(x) }
&=&\sum_{k=0}^n(-1)^{k-1}(x+1)\ldots(x+q^{k-1})q^k
           \frac{(z;q)_k}{(q;q)_k}\frac{(zq^{k+1};q)_{n-k}}{(q;q)_{n-k}}\\
&=&\frac{(z;q)_{n+1}}{(q;q)_n}\sum_{k=0}^n{n \brack
                  k}\frac{(-1)^{k-1}(x+1)\ldots(x+q^{k-1})}{1-zq^k}q^k,
\end{eqnarray*}
and this proves (\ref{EQ2}). \qed

\section{Special Cases}

In their study of overpartitions \cite[Theorem 4.4]{CorLove},
Corteel and Lovejoy obtained a combinatorial interpretation of the
identity:
\begin{eqnarray}    \label{CortLovejoy}
\sum_{i=1}^{\infty}(-1)^{i-1}\frac{(-1;q)_i}{(q;q)_i}\frac{q^i}{1-q^i}
=\sum_{i=1}^{\infty}\frac{2q^{2i-1}}{1-q^{2i-1}}
=\sum_{i=1}^{\infty}\frac{2q^i}{1-q^{2i}} \,  .
\end{eqnarray}

Corteel asked us a proof of the following two related identities~:
\begin{multline}\label{Conje1}
\sum_{i=1}^{2n-1}{2n-1 \brack i}\frac{(-1)^{i-1}
                                       (-1;q)_i}{(1-q^i)^m}q^{mi}\\
=\sum_{i=1}^n\frac{2q^{2i-1}}{1-q^{2i-1}}
      \sum_{2i-1\leq i_2\leq\ldots\leq i_m\leq 2n-1}
          \frac{q^{\sum_{j=2}^mi_j}}{\prod_{j=2}^m(1-q^{i_j})},
\end{multline}
and
\begin{equation}\label{Conje2}
\sum_{i=1}^{2n}{2n \brack i}
           \frac{(-1)^{i-1}(-1;q)_i}{1-q^{i+2}}q^i
=\sum_{i=1}^{n}\frac{2q^{2i-1}(1-q)}{(1+q^2)(1-q^{2i-1})(1-q^{2i+1})}\, ,
\end{equation}
the first one becoming (\ref{CortLovejoy}) for $m=1$ and $n=\infty$,
because
$$   \sum_{i=1}^{\infty}\frac{2q^i}{1-q^{2i}}  =  
  \sum_{i=1}^{\infty} \frac{2q^{2i-1}}{1-q^{2i-1}} \, . $$

In fact, (\ref{Conje1}) is the special case of (\ref{EQ1}) when
$x=1$.  

Changing $n$ into $2n$, specializing $x=1$ and $z=q^2$, one gets from 
(\ref{EQ2}) 
$$ \frac{(1-q^{2n+1})(1-q^{2n+2})}{1-q} 
\sum_0^{2n} {2n \brack i}
           \frac{(-1)^{i-1}(-1;q)_i}{1-q^{i+2}}q^i  =
\sum_0^{2n} (-1)^{i-1} q^i \frac{1-q^{i+1}}{1-q} \, . 
  $$
Suppressing in the left member the term corresponding to $i=0$, 
and taking into account that the right member of (\ref{Conje2}) sums to
$2q (1-q^{2n}) (1-q^4)^{-1} (1-q^{2n+1})^{-1}$, 
one obtains (\ref{Conje2}).

The case $x=0$, $m=1$ of (\ref{EQ1}) is due to Van Hamme
\cite{Hamme} (see also \cite{Andrews}, \cite{FuLa},
\cite{Uchimura}):
\begin{eqnarray*}
\sum_{i=1}^n{n \brack
i}\frac{(-1)^{i-1}q^{\binom{i+1}{2}}}{1-q^i}=\sum_{i=1}^n\frac{q^i}{1-q^i}.
\end{eqnarray*}

Taking $x=0$, and (\ref{Pas}), we get the formula of Dilcher
\cite{Dilcher}:
\begin{eqnarray*}
\sum_{i=1}^n{n \brack i}(-1)^{i-1}\frac{q^{{i \choose
2}+mi}}{(1-q^i)^m}= \sum_{1\leq i_1\leq i_2\leq \ldots \leq i_m
\leq n}\frac{q^{i_1}}{1-q^{i_1}}\ldots\frac{q^{i_m}}{1-q^{i_m}}.
\end{eqnarray*}

When $x=0$ and $z=q^m$ in (\ref{EQ2}), we get Uchimura's identity
\cite{Uchimura2}:
\begin{eqnarray*}
\sum_{i=1}^n{n \brack
i}\frac{(-1)^{i-1}q^{\binom{i+2}{2}}}{1-q^{i+m}}=\sum_{i=1}^n
\frac{q^i}{1-q^i}\left/ {i+m \brack i}.\right.
\end{eqnarray*}

\medskip
\centerline{Acknowledgements}
We thank Sylvie Corteel and Jeremy Lovejoy for their comments.

\end{document}